\documentclass{article}
\usepackage{amsfonts}

\newtheorem{theorem}{Theorem}[section]

\newtheorem{corollary}[theorem]{Corollary}

\newtheorem{lemma}[theorem]{Lemma}

\usepackage{mathrsfs}
\usepackage{amssymb}
\begin{document}

\title{Gr\"obner-Shirshov basis for the  finitely presented  algebras
defined by permutation relations of symmetric
type\footnote{Supported by the NNSF of China (11171118), the
Research Fund for the Doctoral Program of Higher Education of China
(20114407110007), the NSF of Guangdong Province (S2011010003374,S2012040007369),
the Program on International Cooperation and Innovation, Department
of Education, Guangdong Province (2012gjhz0007), and the NSF of
Zhanjiang Normal University (QL0902).}}

\author{
Jianjun Qiu\\
{\small \ Mathematics and Computational  Science School,}
{\small Zhanjiang  Normal University}\\
{\small \ Zhanjiang  524048, China}\\
{\small \ jianjunqiu@126.com}\\
 Yuqun
Chen\footnote {Corresponding author.}\\
{\small \ School of Mathematical Sciences, South China Normal
University}\\
{\small Guangzhou 510631, China}\\
{\small yqchen@scnu.edu.cn}\\
\\
} \vspace{4mm}

\date{}
 \maketitle

\begin{abstract}\noindent
In this paper,  we give  a Gr\"obner-Shirshov basis for the finitely
presented semigroup algebra $\mathbf{k}[S_n(Sym_n)]$ defined by
permutation relations of symmetric type.  As an application, by  the
Composition-Diamond Lemma, we  obtain normal forms  of  elements of
momoid $S_n(Sym_n)$, which gives an answer to an open problem posted
by  F. Ced\'{o}, E. Jespers and J. Okni\'{n}ski \cite{cjo10}  for
the  symmetric group  case.
\vspace{4mm}
 \vspace{4mm}
\end{abstract}

\section{Introduction}

Let $Sym_n$ be the symmetric group of degree $n$ and  $H$  a subset
of $Sym_n$. Recently, F. Ced\'{o}, E. Jespers and J. Okni\'{n}ski
\cite{cjo10} introduced a new class of finitely presented  semigroup
algebra  $\mathbf{k}[S_n(H)]$  over a field $\mathbf{k}$, where the
monoid $S_n(H)$ is defined by  a set of generators
$x_1,x_2,\ldots,x_n$ and   homogenous permutation relations, i.e.
$$
S_n(H)=\langle x_1,x_2,\ldots,x_n|x_{\sigma(1)}x_{\sigma(2)}
\cdots x_{\sigma(n)}=x_1x_2\cdots x_n, \sigma\in H \rangle.
$$
There are some results on this  new algebraic structure,  for
example, the alternating type \cite{cjo09, cjo101},  the abelian
type \cite{cjo12}, and the $n$-cyclic type \cite{cjo10}.

Let $\varsigma$ be the cyclic permutation
\begin{eqnarray}\label{eq1}
\varsigma=\left(
  \begin{array}{ccccc}
    1 & 2 & \cdots & n-1 & n \\
    2 & 3 & \cdots & n& 1 \\
  \end{array}
\right)
\end{eqnarray}
By using the rewriting system method,   Ced\'{o},  Jespers and
Okni\'{n}ski \cite{cjo10} obtained normal forms of elements of
$S_n(H)$  for the case when $H$ is the cyclic subgroup of $Sym_n$
generated by the cyclic permutation $\varsigma$. They also proposed
some open problems at the end of the same   paper \cite{cjo10}. One
of the open problems  is: ``For an arbitrary subgroup $H$ of
symmetric group $Sym_n$,  what does every element of $S_{n}(H)$ have
a unique canonical form, as is the case of the monoid  defined by
permutation relations of cyclic subgroup type."

In this  paper, we  use the  Gr\"obner-Shirshov bases method to
study the finitely presented algebra defined by permutation
relations of symmetric type  $\mathbf{k}[S_n(Sym_n)]$.  We find a
Gr\"obner-Shirshov basis for the algebra $\mathbf{k}[S_n(Sym_n)]$.
As an application, we get normal forms of  elements of monoid
$S_n(Sym_n)$, which gives an answer to the above problem for the
case when $H$ is the symmetric group $Sym_n$.

\section{Composition-Diamond Lemma for associative algebra}

We first cite some concepts and results from the literature
\cite{b72,Sh} which are related to Gr\"{o}bner-Shirshov bases for
associative algebras.

Let $\mathbf{k}$ be a field, $\mathbf{k}\langle X\rangle$ the free associative algebra
over $\mathbf{k}$ generated by $X$. Denote $X^*$ the free monoid generated by
$X$, where the empty word is the identity which is denoted by 1. For
a word $w\in X^*$, we denote the length of $w$ by $|w|$. Let $X^*$
be a well ordered set. Then every nonzero  polynomial  $f\in
\mathbf{k}\langle X\rangle$ has the leading word $\bar{f}$. If the
coefficient of $\bar{f}$ in $f$ is equal to 1, then  $f$ is called
monic.

Let $f$ and $g$ be two monic polynomials in $\mathbf{k}\langle X\rangle$.
Then, there are two kinds of compositions:

$(i)$ If \ $w$ is a word such that $w=\bar{f}b=a\bar{g}$ for some
$a,b\in X^*$ with $|\bar{f}|+|\bar{g}|>|w|$, then the polynomial
 $(f,g)_w=fb-ag$ is called the intersection composition of $f$ and
$g$ with respect to $w$.

$(ii)$ If  $w=\bar{f}=a\bar{g}b$ for some $a,b\in X^*$, then the
polynomial $(f,g)_w=f - agb$ is called the inclusion composition of
$f$ and $g$ with respect to $w$.

In (i) and (ii), the word $w$ is called an ambiguity.

Let $S\subseteq$ $\mathbf{k}\langle X\rangle$ with each $s\in S$ monic. Then
the composition $(f,g)_w$ is called trivial modulo $(S,\ w)$ if
$(f,g)_w=\sum\alpha_i a_i s_i b_i$, where each $\alpha_i\in \mathbf{k}$,
$a_i,b_i\in X^{*}, \ s_i\in S$ and $a_i \overline{ s_i }b_i<w$. If
this is the case, then we write
$$
(f,g)_w\equiv0\quad mod(S,w).
$$
In general, for $p,q\in \mathbf{k}\langle X\rangle$, we write
$$
p\equiv q\quad mod (S,w)
$$
which means that $p-q=\sum\alpha_i a_i s_i b_i $, where each
$\alpha_i\in \mathbf{k},a_i,b_i\in X^{*}, \ s_i\in S$ and $a_i \overline{
s_i} b_i<w$.

We call the set $S$ endowed with the well order $<$ a
Gr\"{o}bner-Shirshov basis  in   $\mathbf{k}\langle X\rangle$ if any
composition of polynomials in $S$ is trivial modulo $S$ and
corresponding $w$.

A well order  $<$ on $X^*$ is monomial if for $u, v\in X^*$, we have
$$
u < v \Rightarrow w_{1}uw_{2} < w_{1}vw_{2},  \ for \  all \
 w_{1}, \ w_{2}\in  X^*.
$$

The following lemma was proved by Shirshov \cite{Sh} for  free Lie
algebras (with deg-lex order) in 1962 (see also Bokut \cite{b72}).
In 1976, Bokut \cite{b76} specialized the approach of Shirshov to
associative algebras (see also Bergman \cite{b}). For commutative
polynomials, this lemma is known as the Buchberger's Theorem (see
\cite{bu65} and \cite{bu70}).

\ \

\noindent{\bf Composition-Diamond Lemma.} \ Let $\mathbf{k}$ be a field, $\mathbf{k}
\langle X|S\rangle=\mathbf{k}\langle X\rangle/Id(S)$ and $>$ a monomial order
on $X^*$, where $Id(S)$ is the ideal of $\mathbf{k} \langle X\rangle$
generated by $S$. Then the following statements are equivalent:
\begin{enumerate}
\item[(i)] $S $ is a Gr\"{o}bner-Shirshov basis  in $\mathbf{k} \langle X\rangle$.
\item[(ii)] $f\in Id(S)\Rightarrow \bar{f}=a\bar{s}b$
for some $s\in S$ and $a,b\in  X^*$.
\item[(iii)] $Irr(S) = \{ u \in X^* |  u \neq a\bar{s}b ,s\in S,a ,b \in X^*\}$
is a $\mathbf{k}$-linear basis of the algebra $\mathbf{k}\langle X |
S \rangle$.
\end{enumerate}

If a subset $S$ of $\mathbf{k}\langle X\rangle$ is not a Gr\"{o}bner-Shirshov
basis, then we can add to $S$ all nontrivial compositions of
polynomials of $S$, and by continuing this process (may be
infinitely) many times, we eventually obtain a Gr\"{o}bner-Shirshov
basis $S^{comp}$. Such a process is called the Shirshov algorithm.

Let $M=\langle X|S\rangle$ be a monoid  presentation. Then $S$ is a
subset of $\mathbf{k}\langle X\rangle$  and hence one can find a
Gr\"{o}bner-Shirshov basis $S^{comp}$.  We also call $S^{comp}$ a
Gr\"{o}bner-Shirshov basis of monoid $M$. The  set
$Irr(S^{comp})=\{u\in X^*|u\neq a\overline{s}b,\ a ,b \in X^*,\ s\in
S^{comp}\}$ is a $\mathbf{k}$-linear basis of $\mathbf{k}\langle
X|S\rangle$ which is also normal forms of elements of monoid $M$.

\section{A Gr\"{o}bner-Shirshov basis for  $\mathbf{k}[S_n(Sym_n)]$}

Let $S_n(Sym_n)$ be  the  finitely presented  momoid  defined by
permutation relations of symmetric type, i.e.
$$
S_n(Sym_n)=\langle x_1,x_2,\ldots,x_n|x_{\sigma(1)}x_{\sigma(2)}
\cdots x_{\sigma(n)}=x_1x_2\cdots x_n, \sigma\in Sym_n\rangle,
$$
where $Sym_n$ is the symmetric group of degree $n$.

We give some notations which will be used in this section. Let
$\varepsilon\in Sym_n$  be  the identity map of $Sym_n$ and
$Sym_n^0=Sym_n\backslash \{\varepsilon\}$. Let $\mathbb{N}$ be the
set of positive integers. Denote $\mathbf{n}=\{1,2,\ldots,n\}$, and
$[n_1,n_2]=\{n_1, n_1+1,\ldots, n_2\}$ for any $n_1,n_2\in\mathbf{n}
$ and  $n_1\leq n_2$.  For  any  $\sigma \in Sym_n$, denote
$$
\mathbf{x}_{\sigma}:=x_{\sigma(1)}x_{\sigma(2)}\cdots x_{\sigma(n)},
$$
in particular,
$$
 \mathbf{x}_{\varepsilon}:=x_1x_2\cdots x_n.
$$
For any $x_{i_1},\ x_{i_2},\ \cdots, x_{i_m}\in X,\ m\geq2$, define
$$
\underline{x_{i_1}x_{i_2}\cdots x_{i_m}}:=x_{j_1}x_{j_2}\cdots
x_{j_m},
$$
where $j_1,j_2,\cdots,j_m$ is the permutation of
$i_1,i_2,\cdots,i_m$ such that $j_1\leq j_2\leq\cdots\leq j_m$. For
example, $\underline{x_2x_5x_4x_3x_2x_3}=x_2x_2x_3x_3x_4x_5$.

Let  $X=\{ x_1,x_2,\ldots,x_n\}$, $x_1<x_2<\cdots<x_n$ and $``<"$
the degree-lexicographic order on $X^*$. Denote
$$
S=\{\mathbf{x}_{\sigma}-\mathbf{x}_{\varepsilon}|\sigma \in Sym_n^0\}
$$
and   $\widetilde{S}$ the subset of $\mathbf{k}\langle X\rangle$
 consisting  of  the following polynomials:
\begin{enumerate}
\item[1]\ $\mathbf{x}_{\sigma}-\mathbf{x}_{\varepsilon},$
\item[2]\
$x_i\mathbf{x}_{\varepsilon}-\mathbf{x}_{\varepsilon}x_i$,
  \item[3]\
 $x_ix_1^m\mathbf{x}_{\varepsilon}-x_1^m\mathbf{x}_{\varepsilon}x_i$,
  \item[4]\
 $\mathbf{x}_{\varepsilon} x_{i_1}x_{i_2}\cdots x_{i_{m+1}}-
 \mathbf{x}_{\varepsilon}\underline{x_{i_1}x_{i_2}\cdots x_{i_{m+1}}}$,
 $x_{i_1}x_{i_2}\cdots x_{i_{m+1}}>\underline{x_{i_1}x_{i_2}\cdots
 x_{i_{m+1}}}$,
   \item[5]\
 $\mathbf{x}_{\varepsilon}x_{i_1}x_{i_2}\cdots x_{i_m}x_1-x_1
 \mathbf{x}_{\varepsilon}x_{i_1}x_{i_2}\cdots x_{i_m}$,

 \end{enumerate}
where $ \sigma \in Sym_n^0$,\ $m\geq 1,\ 2\leq i, i_1, i_2,\cdots,
i_{m+1} \leq n$.

 \begin{lemma}
 $\mathbf{k}[S_n(Sym_n)]=\mathbf{k}\langle X|S\rangle=\mathbf{k}\langle
 X|\widetilde{S}\rangle$.
\end{lemma}

{\bf Proof.}  For any $s_1,s_2\in \mathbf{k}\langle X\rangle$,  we
write $s_1\equiv_I s_2$  if $s_1-s_2\in Id(S)$.  Since $S\subseteq
\widetilde{S}$, we just have to prove that $\widetilde{S}\subseteq
Id(S)$. It suffices to prove that $s\equiv_I 0$ for any $s\in
\widetilde{S}$.

For $2\leq i\leq n$, there exist $\sigma_1,  \sigma_2\in Sym_n^0$
such that $ \mathbf{x}_{\sigma_1}x_i=x_i \mathbf{x}_{\sigma_2}$.
Therefore
$$
x_i\mathbf{x}_{\varepsilon}-\mathbf{x}_{\varepsilon}x_i=
(\mathbf{x}_{\sigma_1}-\mathbf{x}_{\varepsilon})x_i-
x_i(\mathbf{x}_{\sigma_2}-\mathbf{x}_{\varepsilon})\equiv_I0.
$$

Now we use induction on $m$ to prove that all the polynomials of type 3, 4, 5 are in $Id(S)$.

(a)  For $m=1$ and $2\leq i\leq n$, there exist $\sigma_1,
\sigma_2\in Sym_n^0$  such that
$\mathbf{x}_{\sigma_1}x_1x_i=x_ix_1\mathbf{x}_{\sigma_2}$. Therefore
\begin{eqnarray*}
&& x_ix_1\mathbf{x}_{\varepsilon}-x_1\mathbf{x}_{\varepsilon}x_i\\
&=&(\mathbf{x}_{\sigma_1}-\mathbf{x}_{\varepsilon})x_1x_i-x_ix_1(\mathbf{x}_{\sigma_2}-\mathbf{x}_{\varepsilon})+x_1(\mathbf{x}_{\varsigma}-\mathbf{x}_{\varepsilon})x_i\\
&\equiv_I& 0,
\end{eqnarray*}
where $\varsigma$ is  the  cyclic permutation defined by (\ref{eq1}).

(b)  For $m=1$ and $2\leq i_2 < i_1 \leq n$, there exist $\sigma_1,
\sigma_2\in Sym_n^0$ such that
$\mathbf{x}_{\sigma_1}x_{i_2}x_{i_1}=x_{i_1}x_{i_2}\mathbf{x}_{\sigma_2}.$
Therefore,
\begin{eqnarray*}
&&\mathbf{x}_{\varepsilon}x_{i_1}x_{i_2}-\mathbf{x}_{\varepsilon}\underline{x_{i_1}x_{i_2}}\\
&=&\mathbf{x}_{\varepsilon}x_{i_1}x_{i_2}-\mathbf{x}_{\varepsilon}x_{i_2}x_{i_1}\\
&\equiv_I& x_{i_1}x_{i_2}\mathbf{x}_{\varepsilon}-\mathbf{x}_{\varepsilon}x_{i_2}x_{i_1} \ (\mbox{by type  2})\\
&\equiv_I&(\mathbf{x}_{\sigma_1}-\mathbf{x}_{\varepsilon})x_{i_2}x_{i_1}-x_{i_1}x_{i_2}(\mathbf{x}_{\sigma_2}-\mathbf{x}_{\varepsilon})\\
&\equiv_I& 0.
\end{eqnarray*}

(c)  For $m=1$ and  $2\leq i_1\leq n$, since
$\mathbf{x}_{\varepsilon}x_1\equiv_I x_1\mathbf{x}_{\varepsilon}$,
we have
\begin{eqnarray*}
&&\mathbf{x}_{\varepsilon}x_{i_1}x_1-x_1\mathbf{x}_{\varepsilon}x_{i_1}\\
&=& (\mathbf{x}_{\varepsilon}x_{i_1}-x_{i_1}\mathbf{x}_{\varepsilon})x_1+x_{i_1}
\mathbf{x}_{\varepsilon}x_1-x_1\mathbf{x}_{\varepsilon}x_{i_1}\\
&\equiv_I&
(\mathbf{x}_{\varepsilon}x_{i_1}-x_{i_1}\mathbf{x}_{\varepsilon})x_1+(x_{i_1}x_1
\mathbf{x}_{\varepsilon}-x_1\mathbf{x}_{\varepsilon}x_{i_1})\\
&\equiv_I&0  \ (\mbox{by (a) and type 2}).
\end{eqnarray*}

Now we assume that all the polynomials of type 3, 4, 5 are in
$Id(S)$ for $m',\ 1\leq m'< m$.

(i) For $2\leq i\leq n$,  since
$x_1\mathbf{x}_{\varsigma}=\mathbf{x}_{\varepsilon}x_1$, we have
\begin{eqnarray*}
&&x_ix_1^m\mathbf{x}_{\varepsilon}-x_1^m\mathbf{x}_{\varepsilon}x_i\\
&=& x_ix_1^m(-\mathbf{x}_{\varsigma}+\mathbf{x}_{\varepsilon})+ x_ix_1^m\mathbf{x}_{\varsigma}-x_1^m\mathbf{x}_{\varepsilon}x_i\\
&\equiv_I & x_ix_1^{m-1}\mathbf{x}_{\varepsilon} x_1-x_1^m\mathbf{x}_{\varepsilon}x_i \\
&\equiv_I &x_1^{m-1}\mathbf{x}_{\varepsilon}x_ix_1-x_1^m\mathbf{x}_{\varepsilon}x_i \ (\mbox{by induction})\\
&\equiv_I& x_1^{m-1}x_1\mathbf{x}_{\varepsilon}x_i-x_1^m\mathbf{x}_{\varepsilon}x_i\\
&\equiv_I& 0.
\end{eqnarray*}

This shows that all polynomials of type 3 are in $Id(S)$.

(ii) For $2\leq i_1, i_2, \cdots, i_{m+1}\leq n$,    let
$$
x_{i_t}=\max\{x_{i_1},x_{i_2},\cdots, x_{i_{m+1}}\}.
$$
There are two cases to consider.

Case 1: If $x_{i_t}=x_{i_{m+1}}$, then
\begin{eqnarray*}
&& \mathbf{x}_{\varepsilon}x_{i_1}x_{i_2}\cdots x_{i_{m}} x_{i_{m+1}}-
\mathbf{x}_{\varepsilon}\underline{x_{i_1}x_{i_2}\cdots x_{i_{m}} x_{i_{m+1}}}\\
&\equiv_I &\mathbf{x}_{\varepsilon}\underline{x_{i_1}x_{i_2}\cdots x_{i_{m}}}x_{i_{m+1}}-
\mathbf{x}_{\varepsilon}\underline{x_{i_1}x_{i_2}\cdots x_{i_{m}} x_{i_{m+1}}} \ (\mbox{by induction})\\
&\equiv_I &\mathbf{x}_{\varepsilon}\underline{x_{i_1}x_{i_2}
\cdots x_{i_{m}}x_{i_{m+1}}}-\mathbf{x}_{\varepsilon}\underline{x_{i_1}x_{i_2}\cdots x_{i_{m}} x_{i_{m+1}}}\\
&\equiv_I & 0.
\end{eqnarray*}

Case 2: If  $x_{it}\neq x_{i_{m+1}}$, then by induction, we have
\begin{eqnarray*}
&& \mathbf{x}_{\varepsilon}x_{i_1}\cdots x_{i_{t}}\cdots x_{i_{m+1}}-
\mathbf{x}_{\varepsilon}\underline{x_{i_1}\cdots x_{i_{t}}\cdots x_{i_{m+1}}}\\
&\equiv_I & x_{i_1}\cdots x_{i_{t-1}}\mathbf{x}_{\varepsilon}x_{i_{t}}\cdots x_{i_{m+1}}-
\mathbf{x}_{\varepsilon}\underline{x_{i_1}\cdots x_{i_{t}}\cdots x_{i_{m+1}}} \ (\mbox{by\ type\ 2 })\\
&\equiv_I &x_{i_1}\cdots x_{i_{t-1}}\mathbf{x}_{\varepsilon}\underline{x_{i_{t}}\cdots x_{i_{m+1}}}-
\mathbf{x}_{\varepsilon}\underline{x_{i_1}\cdots x_{i_{t}}\cdots x_{i_{m+1}}} \ (\mbox{by induction})\\
&\equiv_I &x_{i_1}\cdots x_{i_{t-1}}\mathbf{x}_{\varepsilon}\underline{x_{i_{t+1}}\cdots x_{i_{m+1}}}x_{i_{t}}-   \mathbf{x}_{\varepsilon}\underline{x_{i_1}\cdots x_{i_{t}}\cdots x_{i_{m+1}}}\\
&\equiv_I &\mathbf{x}_{\varepsilon}x_{i_1}\cdots x_{i_{t-1}}\underline{x_{i_{t+1}}\cdots x_{i_{m+1}}}x_{i_{t}}-  \mathbf{x}_{\varepsilon}\underline{x_{i_1}\cdots x_{i_{t}}\cdots x_{i_{m+1}}}\\
&\equiv_I&\mathbf{x}_{\varepsilon}\underline{x_{i_1}\cdots
x_{i_{t-1}} \underline{x_{i_{t+1}}\cdots x_{i_{m+1}}}}x_{i_{t}}-
\mathbf{x}_{\varepsilon}
\underline{x_{i_1}\cdots x_{i_{t}}\cdots x_{i_{m+1}}}\\
&\equiv_I& 0.
\end{eqnarray*}

This shows that all polynomials of type 4 are in $Id(S)$.

(iii) For  $2\leq i_1, i_2, \cdots, i_m\leq n$, we have,
\begin{eqnarray*}
&&\mathbf{x}_{\varepsilon}x_{i_1}x_{i_2}\cdots x_{i_m}x_1-
x_1\mathbf{x}_{\varepsilon}x_{i_1}x_{i_2}\cdots x_{i_m}\\
&\equiv_I& x_{i_1}\mathbf{x}_{\varepsilon}x_{i_2}\cdots x_{i_m}x_1-
x_1\mathbf{x}_{\varepsilon}x_{i_1}x_{i_2}\cdots x_{i_m} \ (\mbox{by type 2 })\\
&\equiv_I& x_{i_1}x_1\mathbf{x}_{\varepsilon}x_{i_2}\cdots x_{i_m}-
x_1\mathbf{x}_{\varepsilon}x_{i_1}x_{i_2}\cdots x_{i_m} \ (\mbox{by induction})\\
&\equiv_I& x_1 \mathbf{x}_{\varepsilon}x_{i_1}x_{i_2}\cdots x_{i_m}-
x_1\mathbf{x}_{\varepsilon}x_{i_1}x_{i_2}\cdots x_{i_m} \ (\mbox{by type 3 })\\
&\equiv_I & 0.
\end{eqnarray*}

This shows that all polynomials of type 5 are in $Id(S)$.

The proof is complete. \hfill $\blacksquare$

\ \

The following theorem is the main result in this paper.

\begin{theorem}\label{th3.2}
With the degree-lexicographic order on $X^*$, $\widetilde{S}$ is a
Gr\"obner-Shirshov basis in $\mathbf{k}\langle X\rangle$.
\end{theorem}

{\bf Proof. } Let $f_i$ or $f_i'$ be the polynomial of type $i$ in
$\widetilde{S}$, $i=1,2,\ldots, 5$ and $\sigma, \sigma'\in Sym_n^0$.

Denote  $i\wedge j$ the composition of the polynomials of type $i$
and type $j$.

All  possible compositions of the polynomials in $\widetilde{S}$ are
only as below:

\ \

$1\wedge 1$, $f_1=\mathbf{x}_{\sigma}-\mathbf{x}_{\varepsilon},
 f_1'=\mathbf{x}_{\sigma'}-\mathbf{x}_{\varepsilon},$
 $w=x_{i_1}x_{i_2}\cdots x_{i_r}\Delta x_{i_{\pi(1)}}x_{i_{\pi(2)}}\cdots x_{i_{\pi(r)}},$
 $\mathbf{x}_{\sigma}=x_{i_1}x_{i_2}\cdots x_{i_r}\Delta $,
 $\mathbf{x}_{\sigma'}=\Delta x_{i_{\pi(1)}}x_{i_{\pi(2)}}\cdots x_{i_{\pi(r)}}$,
   $\pi\in Sym_r$, $1\leq r< n$.\\

$1\wedge 2$,  $f_1=\mathbf{x}_{\sigma}-\mathbf{x}_{\varepsilon}$,
 $f_2=x_i\mathbf{x}_{\varepsilon}-\mathbf{x}_{\varepsilon}x_i$,
$w=x_{i_1}\cdots x_{i_t}x_i\mathbf{x}_{\varepsilon}$,
$\mathbf{x}_{\sigma}=x_{i_1}\cdots x_{i_t} x_i\\x_1x_2\ldots
x_{n-t-1}$, $\{ i_1,i_2,\ldots, i_t\}=[n-t, n]\backslash \{i\}$,
$0\leq n-t-1<i$, $ 2\leq i\leq n$.\\

$1\wedge 3$, there are two cases. Let
$f_1=\mathbf{x}_{\sigma}-\mathbf{x}_{\varepsilon}$ and
$f_3=x_ix_1^m\mathbf{x}_{\varepsilon}-x_1^m\mathbf{x}_{\varepsilon}x_i$.

$w_1=\mathbf{x}_{\sigma}x_1^{m-1}\mathbf{x}_{\varepsilon}$,
 $\mathbf{x}_{\sigma}=x_{i_1}x_{i_2}\cdots x_{i_{n-2}}x_ix_1$,
  $\{i_1, i_2, \ldots, i_{n-2}\}=\mathbf{n}\backslash  \{i, 1\}$, $ m\geq 1$,  $2\leq i\leq n$.

$w_2=\mathbf{x}_{\sigma}x_1^{m}\mathbf{x}_{\varepsilon}$, $\mathbf{x}_{\sigma}=x_{i_1}x_{i_2}\cdots x_{i_{n-1}}x_i$,
$\{i_1, i_2, \ldots, i_{n-1}\}=\mathbf{n}\backslash  \{i\}$, $ m\geq 1$,  $2\leq i\leq n$.\\

$1\wedge 4$, $f_1=\mathbf{x}_{\sigma}-\mathbf{x}_{\varepsilon}$,
$f_4=\mathbf{x}_{\varepsilon}x_{i_1}x_{i_2}\cdots
x_{i_{m+1}}-\mathbf{x}_{\varepsilon}\underline{x_{i_1}x_{i_2}\cdots
x_{i_{m+1}}}$, $w=x_{j_1}x_{j_2}\cdots
x_{j_t}\mathbf{x}_{\varepsilon}x_{i_1}x_{i_2}\cdots x_{i_{m+1}}$,
$\mathbf{x}_{\sigma}=x_{j_1}x_{j_2}\cdots x_{j_t}x_1x_2\cdots
x_{n-t}$, $\{j_1, j_2,\ldots, j_t\}=[n-t+1, n]$, $2\leq i_1, i_2,
\cdots,i_{m+1}\leq n,$ $m \geq 1$, $1\leq t\leq n-1$,
$x_{i_1}x_{i_2}\cdots x_{i_{m+1}}> \underline{x_{i_1}x_{i_2}\cdots
x_{i_{m+1}}}$.\\

$1\wedge 5$, $f_1=\mathbf{x}_{\sigma}-\mathbf{x}_{\varepsilon}$,
$f_5=\mathbf{x}_{\varepsilon}x_{i_1}x_{i_2}\cdots
x_{i_m}x_1-x_1\mathbf{x}_{\varepsilon}x_{i_1}x_{i_2}\cdots x_{i_m}$,
$w=x_{j_1}x_{j_2}\cdots
x_{j_t}\mathbf{x}_{\varepsilon}x_{i_1}x_{i_2}\cdots x_{i_{m}}x_1$,
$\mathbf{x}_{\sigma}=x_{j_1}x_{j_2}\cdots x_{j_t}x_1x_2\cdots
x_{n-t}$,
$\{j_1, j_2,\ldots, j_t\}=[n-t+1, n]$, $2\leq i_1, i_2, \cdots, i_m\leq n$, $m \geq 1$, $1\leq t\leq n-1$.\\

$2\wedge 1$,
$f_2=x_i\mathbf{x}_{\varepsilon}-\mathbf{x}_{\varepsilon}x_i,$
$f_1=\mathbf{x}_{\sigma}-\mathbf{x}_{\varepsilon}$,
$w=x_i\mathbf{x}_{\varepsilon}x_{j_1}x_{j_2}\cdots x_{j_t}$,
 $\mathbf{x}_{\sigma}=x_{t+1}\cdots x_n x_{j_1}x_{j_2}\cdots x_{j_t}$,
 $\{j_1,j_2, \ldots, j_t \}=\mathbf{n}\backslash [t+1, n]$, $2\leq i\leq n$,  $1\leq t\leq n-1$.\\

$2\wedge 2$,
$f_2=x_i\mathbf{x}_{\varepsilon}-\mathbf{x}_{\varepsilon}x_i,$,
$f_2'=x_n\mathbf{x}_{\varepsilon}-\mathbf{x}_{\varepsilon}x_n$,
$w=x_i\mathbf{x}_{\varepsilon}\mathbf{x}_{\varepsilon}$, $2\leq
i\leq  n$.\\

$2\wedge 3$,
$f_2=x_i\mathbf{x}_{\varepsilon}-\mathbf{x}_{\varepsilon}x_i$,
 $f_3=x_nx_1^m\mathbf{x}_{\varepsilon}-x_1^m\mathbf{x}_{\varepsilon}x_n$,
  $w=x_i\mathbf{x}_{\varepsilon}x_1^m\mathbf{x}_{\varepsilon}$,  $2 \leq i\leq n$, $m\geq 1$.\\

$2\wedge 4$,
$f_2=x_i\mathbf{x}_{\varepsilon}-\mathbf{x}_{\varepsilon}x_i$,
 $f_4=\mathbf{x}_{\varepsilon}x_{i_1}x_{i_2}\cdots x_{i_{m+1}}-\mathbf{x}_{\varepsilon}\underline{x_{i_1}x_{i_2}\cdots x_{i_{m+1}}}$,
$w=x_i\mathbf{x}_{\varepsilon}x_{i_1}x_{i_2}\cdots x_{i_{m+1}}$,
 $m\geq 1$, $2\leq i, i_1,i_2,\ldots,i_{m+1}\leq n$,
   $x_{i_1}x_{i_2}\cdots x_{i_{m+1}}> \underline{x_{i_1}x_{i_2}\cdots
x_{i_{m+1}}}$.\\

$2\wedge 5$,
$f_2=x_i\mathbf{x}_{\varepsilon}-\mathbf{x}_{\varepsilon}x_i$,
$f_5=\mathbf{x}_{\varepsilon}x_{i_1}x_{i_2}\cdots
x_{i_m}x_1-x_1\mathbf{x}_{\varepsilon}x_{i_1}x_{i_2}\cdots x_{i_m}$,
 $w=x_i\mathbf{x}_{\varepsilon}x_{i_1}x_{i_2}\cdots x_{i_m}x_1$, $m\geq 1$,
 $2\leq i, i_1, i_2, \ldots, i_{m}\leq n$.\\

$3\wedge 1$,
$f_3=x_ix_1^m\mathbf{x}_{\varepsilon}-x_1^m\mathbf{x}_{\varepsilon}x_i$,
 $f_1=\mathbf{x}_{\sigma}-\mathbf{x}_{\varepsilon}$, $w=x_ix_1^m\mathbf{x}_{\varepsilon}x_{i_1}\cdots x_{i_t}$,
$\mathbf{x}_{\sigma}=x_{t+1}x_{t+2}\cdots x_{n}x_{i_1}\cdots x_{i_t}$,
 $\{i_1,i_2,\ldots, i_t\}=[1, t]$, $1\leq t\leq n-1$, $2\leq i \leq n$, $m\geq 1$. \\

$3\wedge 2$,
$f_3=x_ix_1^m\mathbf{x}_{\varepsilon}-x_1^m\mathbf{x}_{\varepsilon}x_i$,
$f_2=x_n\mathbf{x}_{\varepsilon}-\mathbf{x}_{\varepsilon}x_n$,
$w=x_ix_1^m\mathbf{x}_{\varepsilon}\mathbf{x}_{\varepsilon}$, $2\leq
i \leq n$, $m\geq 1$.\\

$3\wedge 3$,
$f_3=x_ix_1^m\mathbf{x}_{\varepsilon}-x_1^m\mathbf{x}_{\varepsilon}x_i$,
  $f_3'=x_nx_1^{m_1}\mathbf{x}_{\varepsilon}-x_1^{m_1}\mathbf{x}_{\varepsilon}x_n$,
  $w=x_ix_1^{m}\mathbf{x}_{\varepsilon}x_1^{m_1}\mathbf{x}_{\varepsilon}$,
$2\leq i \leq n$, $m, m_1\geq 1$.\\

$3\wedge 4$,
$f_3=x_ix_1^m\mathbf{x}_{\varepsilon}-x_1^m\mathbf{x}_{\varepsilon}x_i$,
$f_4=\mathbf{x}_{\varepsilon}x_{i_1}\cdots
x_{i_{m_1+1}}-\mathbf{x}_{\varepsilon}\underline{x_{i_1}\cdots
x_{i_{m_1+1}}}$,
$w=x_ix_1^m\mathbf{x}_{\varepsilon}x_{i_1}x_{i_2}\cdots
x_{i_{m_1+1}}$,
 $2\leq i, i_1,i_2,\ldots,i_{m_1+1}\leq n$,\
 $m,m_1\geq 1$, $x_{i_1}\cdots
x_{i_{m_1+1}}>\underline{x_{i_1}\cdots
x_{i_{m_1+1}}}$.\\

$3\wedge 5$,
$f_3=x_ix_1^m\mathbf{x}_{\varepsilon}-x_1^m\mathbf{x}_{\varepsilon}x_i$,
 $f_5=\mathbf{x}_{\varepsilon}x_{i_1}x_{i_2}\cdots x_{i_{m_1}}x_1-x_1\mathbf{x}_{\varepsilon}x_{i_1}x_{i_2}\cdots x_{i_{m_1}}$,
$w=x_ix_1^m\mathbf{x}_{\varepsilon}x_{i_1}x_{i_2}\cdots x_{i_{m_1}}x_1$,
 $2\leq i, i_1,i_2,\ldots,i_{m_1}\leq n$,   $m,m_1\geq 1$.\\

$4\wedge 1$, there are two cases. Let
$f_4=\mathbf{x}_{\varepsilon}x_{i_1}x_{i_2}\cdots
x_{i_{m+1}}-\mathbf{x}_{\varepsilon}\underline{x_{i_1}x_{i_2}\cdots
x_{i_{m+1}}}$, $f_1=\mathbf{x}_{\sigma}-\mathbf{x}_{\varepsilon}$,
  $m\geq 1$, $2\leq i_1, i_2,\cdots,i_{m+1}\leq n$,
  $x_{i_1}x_{i_2}\cdots x_{i_{m+1}}>\underline{x_{i_1}x_{i_2}\cdots
x_{i_{m+1}}}$.

$w_1=\mathbf{x}_{\varepsilon}x_{i_1}\cdots
x_{i_{m+1}}x_{j_1}x_{j_2}\cdots x_{j_{t}}$,
$\mathbf{x}_{\sigma}=x_{i_{m+2-n+t}} \cdots
x_{i_{m+1}}x_{j_1}x_{j_2}\cdots x_{j_{t}}$,\\
 $\{ j_1, j_2, \ldots,
j_{t}\}=\mathbf{n}\backslash \{ i_{m+2-n+t}, \ldots, i_{m+1}\}$,

$w_2=\mathbf{x}_{\varepsilon}x_{i_1} \cdots x_{i_{m+1}}x_{j_1}\cdots
x_{j_{t-m-1}}$, $\mathbf{x}_{\sigma}=x_{t+1}\cdots x_n x_{i_1}
\cdots x_{i_{m+1}}x_{j_1}\cdots x_{j_{t-m-1}}$, $2\leq i_1,
i_2,\cdots,i_{m+1}\leq n$, $\{ j_1, j_2, \ldots,
j_{t-m-1}\}=\mathbf{n}\backslash ([t+1, n]\cup \{ i_{1},
i_{2},\ldots, i_{m+1}\})$,
 $t-m-1\geq 1$.\\

$4\wedge 2$, $f_4=\mathbf{x}_{\varepsilon}x_{i_1}x_{i_2}\cdots
x_{i_{m+1}}-\mathbf{x}_{\varepsilon}\underline{x_{i_1}x_{i_2}\cdots
x_{i_{m+1}}}$,
$f_2=x_{i_{m+1}}\mathbf{x}_{\varepsilon}-\mathbf{x}_{\varepsilon}x_{i_{m+1}}$,
$w=\mathbf{x}_{\varepsilon}x_{i_1}x_{i_2}\cdots
x_{i_{m+1}}\mathbf{x}_{\varepsilon}$, $m\geq 1$, $2\leq
i_1,\cdots,i_{m+1}\leq n$, $x_{i_1}x_{i_2}\cdots
x_{i_{m+1}}>\underline{x_{i_1}x_{i_2}\cdots x_{i_{m+1}}}$.\\

$4\wedge 3$, $f_4=\mathbf{x}_{\varepsilon}x_{i_1}x_{i_2}\cdots
x_{i_{m+1}}-\mathbf{x}_{\varepsilon}\underline{x_{i_1}x_{i_2}\cdots
x_{i_{m+1}}}$,
$f_3=x_{i_{m+1}}x_1^{m_1}\mathbf{x}_{\varepsilon}-x_1^{m_1}\mathbf{x}_{\varepsilon}x_{i_{m+1}}$,
$w=\mathbf{x}_{\varepsilon}x_{i_1}x_{i_2}\cdots
x_{i_{m+1}}x_1^{m_1}\mathbf{x}_{\varepsilon}$, $2\leq
i_1,\cdots,i_{m+1}\leq n$, $m_1, m\geq 1$, $x_{i_1}x_{i_2}\cdots
x_{i_{m+1}}>\underline{x_{i_1}x_{i_2}\cdots x_{i_{m+1}}}$.\\

$5\wedge 1$,  there are two cases. Let
 $f_5=\mathbf{x}_{\varepsilon}x_{i_1}x_{i_2}\cdots
x_{i_m}x_1-x_1\mathbf{x}_{\varepsilon}x_{i_1}x_{i_2}\cdots x_{i_m}$,
$f_1=\mathbf{x}_{\sigma}-\mathbf{x}_{\varepsilon}$, $2\leq i_1, i_2,
\cdots , i_m\leq n$, $m\geq 1$.

$w_1=\mathbf{x}_{\varepsilon}x_{i_1}\cdots x_{i_m}x_1x_{j_1}\cdots
x_{j_{n+t-m-2}}$, $\mathbf{x}_{\sigma}=x_{i_{t}}x_{i_{t+1}}\cdots
x_{i_m}x_1x_{j_1}\cdots x_{j_{n+t-m-2}}$, $2\leq i_1, i_2, \cdots ,
i_m\leq n$, $\{ j_1,j_2, \ldots, j_{n+t-m-2} \}
=\mathbf{n}\backslash(\{i_{t},i_{t+1},\ldots,i_m\}\cup \{1\})$.

$w_2=\mathbf{x}_{\varepsilon}x_{i_1}\cdots x_{i_m}x_1$,
$\mathbf{x}_{\sigma}=x_{t+1}\cdots x_nx_{i_1}x_{i_2}\cdots
x_{i_m}x_1x_{j_1}\cdots x_{j_{t-m}}$, \\  $2\leq i_1, i_2, \cdots,
i_m\leq n$, $\{ j_1,j_2, \ldots, j_{t-m}
\}=\mathbf{n}\backslash(\{i_{1},i_{2},\ldots,i_m\}\cup [t+1, n]\cup
\{1\})$, $ 1\leq m \leq n-2$, $t-m\geq 0$.\\

$5\wedge 2$,  $f_5=\mathbf{x}_{\varepsilon}x_{i_1}x_{i_2}\cdots
x_{i_m}x_1-x_1\mathbf{x}_{\varepsilon}x_{i_1}x_{i_2}\cdots x_{i_m}$,
$f_2=x_{i_m}\mathbf{x}_{\varepsilon}-\mathbf{x}_{\varepsilon}x_{i_m}$
$w=\mathbf{x}_{\varepsilon}x_{i_1}x_{i_2}\cdots
x_{i_m}\mathbf{x}_{\varepsilon}$,
$2\leq i_1, i_2, \cdots, i_{m}\leq n$, $m\geq 1$. \\

$5\wedge 3$, $f_5=\mathbf{x}_{\varepsilon}x_{i_1}x_{i_2}\cdots
x_{i_m}x_1-x_1\mathbf{x}_{\varepsilon}x_{i_1}x_{i_2}\cdots x_{i_m}$,
$f_3=
x_{i_m}x_1^{m_1}\mathbf{x}_{\varepsilon}-x_1^{m_1}\mathbf{x}_{\varepsilon}x_{i_m}$,
$w=\mathbf{x}_{\varepsilon}x_{i_1}x_{i_2}\cdots
x_{i_m}x_1^{m_1}\mathbf{x}_{\varepsilon}$,
$2\leq i_1, i_2, \cdots,  i_{m}\leq n$, $m, m_1\geq 1$.\\

$5\wedge 4$,
 $f_5=\mathbf{x}_{\varepsilon}x_{i_1}x_{i_2}\cdots
x_{i_m}x_1-x_1\mathbf{x}_{\varepsilon}x_{i_1}x_{i_2}\cdots x_{i_m}$,
$f_4=\mathbf{x}_{\varepsilon}x_{j_1}x_{j_2}\cdots x_{j_{m_1+1}}$ \\
$-\mathbf{x}_{\varepsilon}\underline{x_{j_1}x_{j_2}\cdots
x_{j_{m_1+1}}}$, $w=\mathbf{x}_{\varepsilon}x_{i_1}x_{i_2}\cdots
x_{i_m}\mathbf{x}_{\varepsilon}x_{j_1}x_{j_2}\cdots x_{j_{m_1+1}}$,
$2\leq i_1, i_2,  \cdots,  i_m,\\ j_1,j_2,\cdots, j_{m_1+1}\leq n,$
 $ m, m_1\geq 1$, $x_{j_1}x_{j_2}\cdots x_{j_{m_1+1}}>\underline{x_{j_1}x_{j_2}\cdots x_{j_{m_1+1}}}$.\\

$5\wedge 5$, $f_5=\mathbf{x}_{\varepsilon}x_{i_1}x_{i_2}\cdots
x_{i_m}x_1-x_1\mathbf{x}_{\varepsilon}x_{i_1}x_{i_2}\cdots x_{i_m}$,
$f_5'=\mathbf{x}_{\varepsilon}x_{i_1'}x_{i_2'}\cdots x_{i_{m_1}'}x_1$ \\
$-x_1\mathbf{x}_{\varepsilon}x_{i_1'}x_{i_2'}\cdots x_{i_{m_1}'}$,
$w=\mathbf{x}_{\varepsilon}x_{i_1}x_{i_2}\cdots x_{i_m}\mathbf{x}_{\varepsilon}x_{i_1'}x_{i_2'}\cdots x_{i_{m_1}'}x_1$,
 $2\leq i_1, i_2,  \cdots,  i_m,\\ i_1', i_2', \cdots,  i_{m_1}'\leq n$.\\

We prove that  all the above  compositions are trivial. Here, we
just check $1\wedge 1$, $1\wedge 4$, $2\wedge 5$. Others are
similarly proved.

For $1\wedge 1$, there are two cases to consider.

Case 1: If $1\notin \{i_1, i_2, \ldots, i_r\}$, then
\begin{eqnarray*}
1\wedge 1&=&f_1x_{i_{\pi(1)}}x_{i_{\pi(2)}}\cdots x_{i_{\pi(r)}}-x_{i_1}x_{i_2}\cdots x_{i_r}f_1'\\
&=&-\mathbf{x}_{\varepsilon}x_{i_{\pi(1)}}x_{i_{\pi(2)}}\cdots x_{i_{\pi(r)}}+
x_{i_1}x_{i_2}\cdots x_{i_r}\mathbf{x}_{\varepsilon}\\
&\equiv &-\mathbf{x}_{\varepsilon}x_{i_{\pi(1)}}x_{i_{\pi(2)}}\cdots x_{i_{\pi(r)}}+
\mathbf{x}_{\varepsilon}x_{i_1}x_{i_2}\cdots x_{i_r} \ (\mbox{by type 2 })\\
&\equiv
&-\mathbf{x}_{\varepsilon}\underline{x_{i_{\pi(1)}}x_{i_{\pi(2)}}\cdots
x_{i_{\pi(r)}}}+
\mathbf{x}_{\varepsilon}\underline{x_{i_1}x_{i_2}\cdots x_{i_r}} \ (\mbox{by type 4 })\\
&\equiv& 0\ mod (\widetilde{S}, w),
\end{eqnarray*}

Case 2: If $1\in \{i_1, i_2, \ldots, i_r\}$, say,
$x_1=x_{i_t}=x_{i_{\pi(s)}}, 1\leq s, t\leq r$, then by type 5 and
4, we have
\begin{eqnarray*}
1\wedge 1&=&f_1x_{i_{\pi(1)}}x_{i_{\pi(2)}}\cdots x_{i_{\pi(r)}}-x_{i_1}x_{i_2}\cdots x_{i_r}f_1'\\
&=&-\mathbf{x}_{\varepsilon}x_{i_{\pi(1)}}\cdots x_{i_{\pi(s)}}\cdots x_{i_{\pi(r)}}+
x_{i_1}\cdots x_{i_t}\cdots x_{i_r}\mathbf{x}_{\varepsilon}\\
&=&-x_1\mathbf{x}_{\varepsilon}x_{i_{\pi(1)}}\cdots x_{i_{\pi(s-1)}} x_{i_{\pi(s+1)}}\cdots x_{i_{\pi(r)}}+
x_1\mathbf{x}_{\varepsilon}x_{i_1}\cdots x_{i_{t-1}}x_{i_{t+1}}\cdots x_{i_r}\\
&=&-x_1\mathbf{x}_{\varepsilon}\underline{x_{i_{\pi(1)}}\cdots x_{i_{\pi(s-1)}} x_{i_{\pi(s+1)}}\cdots x_{i_{\pi(r)}}}+
x_1\mathbf{x}_{\varepsilon}\underline{x_{i_1}\cdots x_{i_{t-1}}x_{i_{t+1}}\cdots x_{i_r}}\\
&\equiv& 0\ mod (\widetilde{S}, w).\\
\\
1\wedge 4&=&f_1x_{n-t+1}x_{n-t+2}\cdots x_nx_{i_1}x_{i_2}\cdots x_{i_{m+1}}-x_{j_1}x_{j_2}\cdots x_{j_t}f_4\\
&=&-\mathbf{x}_{\varepsilon}x_{n-t+1}x_{n-t+2}\cdots x_nx_{i_1}x_{i_2}\cdots x_{i_{m+1}}+
x_{j_1}x_{j_2}\cdots x_{j_t} \mathbf{x}_{\varepsilon}\underline{x_{i_1}x_{i_2}\cdots x_{i_{m+1}}}\\
&\equiv&-\mathbf{x}_{\varepsilon}x_{n-t+1}x_{n-t+2}\cdots x_nx_{i_1}x_{i_2}\cdots x_{i_{m+1}}+
 \mathbf{x}_{\varepsilon}x_{j_1}x_{j_2}\cdots x_{j_t}\underline{x_{i_1}x_{i_2}\cdots x_{i_{m+1}}}\\
&\equiv&-\mathbf{x}_{\varepsilon}\underline{x_{n-t+1}x_{n-t+2}\cdots x_nx_{i_1}x_{i_2}\cdots x_{i_{m+1}}}+
 \mathbf{x}_{\varepsilon}\underline{x_{j_1}x_{j_2}\cdots x_{j_t}\underline{x_{i_1}x_{i_2}\cdots x_{i_{m+1}}}}\\
&\equiv& 0 \ mod(\widetilde{S},w).\\
\\
2\wedge 5&=&f_2x_{i_1}x_{i_2}\cdots x_{i_{m}}x_1-x_if_5\\
&\equiv& \mathbf{x}_{\varepsilon}x_ix_{i_1}x_{i_2}\cdots x_{i_{m}}x_1-x_ix_1\mathbf{x}_{\varepsilon}x_{i_1}x_{i_2}\cdots x_{i_m}\\
&\equiv& \mathbf{x}_{\varepsilon}x_ix_{i_1}x_{i_2}\cdots x_{i_{m}}x_1-x_1 \mathbf{x}_{\varepsilon}x_i x_{i_1}x_{i_2}\cdots x_{i_m}\\
&\equiv& x_1\mathbf{x}_{\varepsilon}x_ix_{i_1}x_{i_2}\cdots x_{i_{m}}-x_1 \mathbf{x}_{\varepsilon}x_i x_{i_1}x_{i_2}\cdots x_{i_m}\\
&\equiv& mod (\widetilde{S}, w).
\end{eqnarray*}

The proof is complete. \hfill $\blacksquare$\\

By the Composition-Diamond Lemma and Theorem \ref{th3.2}, we have the following corollary.
\begin{corollary}
The set
$$
Irr(\widetilde{S})=(X^*\backslash\bigcup_{\sigma \in
Sym_n}X^*\{\mathbf{x}_{\sigma}\}X^*)\bigcup
\{x_1^{m_1}\mathbf{x}_{\varepsilon}x_2^{m_2}\cdots x_n^{m_n}
|m_i\geq 0, i=1,2,\ldots, n\}
$$
is a  $\mathbf{k}$-linear basis of  algebra
$\mathbf{k}[S_n(Sym_n)]$. Moreover, $Irr(\widetilde{S})$ is normal
forms of  elements of monoid $S_n(Sym_n)$.
\end{corollary}

\section{Appendix}
In this section, we will check that all the compositions are trivial.\\

$1\wedge 1$, $f_1=\mathbf{x}_{\sigma}-\mathbf{x}_{\varepsilon},
 f_1'=\mathbf{x}_{\sigma'}-\mathbf{x}_{\varepsilon},$
 $w=x_{i_1}x_{i_2}\cdots x_{i_r}\Delta x_{i_{\pi(1)}}x_{i_{\pi(2)}}\cdots x_{i_{\pi(r)}},$
 $\mathbf{x}_{\sigma}=x_{i_1}x_{i_2}\cdots x_{i_r}\Delta $,
 $\mathbf{x}_{\sigma'}=\Delta x_{i_{\pi(1)}}x_{i_{\pi(2)}}\cdots x_{i_{\pi(r)}}$,
   $\pi\in Sym_r$, $1\leq r< n$.\\

For $1\wedge 1$, there are two cases to consider.

Case 1: If $1\notin \{i_1, i_2, \ldots, i_r\}$, then
\begin{eqnarray*}
1\wedge 1&=&f_1x_{i_{\pi(1)}}x_{i_{\pi(2)}}\cdots x_{i_{\pi(r)}}-x_{i_1}x_{i_2}\cdots x_{i_r}f_1'\\
&=&-\mathbf{x}_{\varepsilon}x_{i_{\pi(1)}}x_{i_{\pi(2)}}\cdots x_{i_{\pi(r)}}+
x_{i_1}x_{i_2}\cdots x_{i_r}\mathbf{x}_{\varepsilon}\\
&\equiv &-\mathbf{x}_{\varepsilon}x_{i_{\pi(1)}}x_{i_{\pi(2)}}\cdots x_{i_{\pi(r)}}+
\mathbf{x}_{\varepsilon}x_{i_1}x_{i_2}\cdots x_{i_r} \ (\mbox{by type 2 })\\
&\equiv
&-\mathbf{x}_{\varepsilon}\underline{x_{i_{\pi(1)}}x_{i_{\pi(2)}}\cdots
x_{i_{\pi(r)}}}+
\mathbf{x}_{\varepsilon}\underline{x_{i_1}x_{i_2}\cdots x_{i_r}} \ (\mbox{by type 4 })\\
&\equiv& 0\ mod (\widetilde{S}, w),
\end{eqnarray*}

Case 2: If $1\in \{i_1, i_2, \ldots, i_r\}$, say,
$x_1=x_{i_t}=x_{i_{\pi(s)}}, 1\leq s, t\leq r$, then by type 5 and
4, we have
\begin{eqnarray*}
1\wedge 1&=&f_1x_{i_{\pi(1)}}x_{i_{\pi(2)}}\cdots x_{i_{\pi(r)}}-x_{i_1}x_{i_2}\cdots x_{i_r}f_1'\\
&=&-\mathbf{x}_{\varepsilon}x_{i_{\pi(1)}}\cdots x_{i_{\pi(s)}}\cdots x_{i_{\pi(r)}}+
x_{i_1}\cdots x_{i_t}\cdots x_{i_r}\mathbf{x}_{\varepsilon}\\
&=&-x_1\mathbf{x}_{\varepsilon}x_{i_{\pi(1)}}\cdots x_{i_{\pi(s-1)}} x_{i_{\pi(s+1)}}\cdots x_{i_{\pi(r)}}+
x_1\mathbf{x}_{\varepsilon}x_{i_1}\cdots x_{i_{t-1}}x_{i_{t+1}}\cdots x_{i_r}\\
&=&-x_1\mathbf{x}_{\varepsilon}\underline{x_{i_{\pi(1)}}\cdots x_{i_{\pi(s-1)}} x_{i_{\pi(s+1)}}\cdots x_{i_{\pi(r)}}}+
x_1\mathbf{x}_{\varepsilon}\underline{x_{i_1}\cdots x_{i_{t-1}}x_{i_{t+1}}\cdots x_{i_r}}\\
&\equiv& 0\ mod (\widetilde{S}, w).
\end{eqnarray*}

$1\wedge 2$,  $f_1=\mathbf{x}_{\sigma}-\mathbf{x}_{\varepsilon}$,
 $f_2=x_i\mathbf{x}_{\varepsilon}-\mathbf{x}_{\varepsilon}x_i$,
$w=x_{i_1}\cdots x_{i_t}x_i\mathbf{x}_{\varepsilon}$,
$\mathbf{x}_{\sigma}=x_{i_1}\cdots x_{i_t} x_i\\x_1x_2\ldots
x_{n-t-1}$, $\{ i_1,i_2,\ldots, i_t\}=[n-t, n]\backslash \{i\}$,
$0\leq n-t-1<i$, $ 2\leq i\leq n$.

\begin{eqnarray*}
1\wedge 2 &=&f_1x_{n-t}\cdots x_{n}-x_{i_1}x_{i_2}\cdots x_{i_t}f_2\\
&=& \mathbf{x}_{\varepsilon}x_{n-t}\cdots x_{n}-x_{i_1}x_{i_2}\cdots x_{i_t}\mathbf{x}_{\varepsilon}x_{i}\\
&=& \mathbf{x}_{\varepsilon}x_{n-t}\cdots x_{n}-\mathbf{x}_{\varepsilon}x_{i_1}x_{i_2}\cdots x_{i_t}x_{i}\\
&\equiv& \mathbf{x}_{\varepsilon}x_{n-t}\cdots x_{n}-\mathbf{x}_{\varepsilon}\underline{x_{i_1}x_{i_2}\cdots x_{i_t}x_{i}}\\
&\equiv& \mathbf{x}_{\varepsilon}x_{n-t}\cdots x_{n}-\mathbf{x}_{\varepsilon}x_{n-t}\cdots x_{n}\\
&\equiv&0\ mod(\widetilde{S}, w)
\end{eqnarray*}

$1\wedge 3$, there are two cases. Let
$f_1=\mathbf{x}_{\sigma}-\mathbf{x}_{\varepsilon}$ and
$f_3=x_ix_1^m\mathbf{x}_{\varepsilon}-x_1^m\mathbf{x}_{\varepsilon}x_i$.

$w_1=\mathbf{x}_{\sigma}x_1^{m-1}\mathbf{x}_{\varepsilon}$,
 $\mathbf{x}_{\sigma}=x_{i_1}x_{i_2}\cdots x_{i_{n-2}}x_ix_1$,
  $\{i_1, i_2, \ldots, i_{n-2}\}=\mathbf{n}\backslash  \{i, 1\}$, $ m\geq 1$,  $2\leq i\leq n$.

$w_2=\mathbf{x}_{\sigma}x_1^{m}\mathbf{x}_{\varepsilon}$, $\mathbf{x}_{\sigma}=x_{i_1}x_{i_2}\cdots x_{i_{n-1}}x_i$,
$\{i_1, i_2, \ldots, i_{n-1}\}=\mathbf{n}\backslash  \{i\}$, $ m\geq 1$,  $2\leq i\leq n$.

\begin{eqnarray*}
1\wedge 3 &=&f_1x_{1}^{m-1}\mathbf{x}_{\varepsilon}-x_{i_1}x_{i_2}\cdots x_{i_{n-2}}f_3\\
&\equiv&\mathbf{x}_{\varepsilon}x_{1}^{m-1}\mathbf{x}_{\varepsilon}-x_{i_1}x_{i_2}\cdots x_{i_{n-2}}x_1^m\mathbf{x}_{\varepsilon}x_i\\
&\equiv&\mathbf{x}_{\varepsilon}x_{1}^{m-1}\mathbf{x}_{\varepsilon}-x_1^m\mathbf{x}_{\varepsilon}x_{i_1}x_{i_2}\cdots x_{i_{n-2}}x_i\\
&\equiv&x_{1}^{m}\mathbf{x}_{\varepsilon}x_2\cdots x_n-x_1^m\mathbf{x}_{\varepsilon}\underline{x_{i_1}x_{i_2}\cdots x_{i_{n-2}}x_i}\\
&\equiv&x_{1}^{m}\mathbf{x}_{\varepsilon}x_2\cdots x_n-x_{1}^{m}\mathbf{x}_{\varepsilon}x_2\cdots x_n\\
&\equiv&0\ mod(\widetilde{S}, w_1)
\end{eqnarray*}

\begin{eqnarray*}
1\wedge 3 &=&f_1x_{1}^{m}\mathbf{x}_{\varepsilon}-x_{i_1}x_{i_2}\cdots x_{i_{n-1}}f_3\\
&\equiv&\mathbf{x}_{\varepsilon}x_{1}^{m}\mathbf{x}_{\varepsilon}-x_{i_1}x_{i_2}\cdots x_{i_{n-1}}x_1^m\mathbf{x}_{\varepsilon}x_i\\
&\equiv&\mathbf{x}_{\varepsilon}x_{1}^{m}\mathbf{x}_{\varepsilon}-x_1^m\mathbf{x}_{\varepsilon}x_{i_1}x_{i_2}\cdots x_{i_{n-1}}x_i\\
&\equiv&x_{1}^{m+1}\mathbf{x}_{\varepsilon}x_2\cdots x_n-x_1^{m+1}\mathbf{x}_{\varepsilon}x_{i_1'}x_{i_2'}\cdots x_{i_{n-2}'}x_i\\
&\equiv&x_{1}^{m+1}\mathbf{x}_{\varepsilon}x_2\cdots x_n-x_{1}^{m+1}\mathbf{x}_{\varepsilon}x_2\cdots x_n\\
&\equiv&0\ mod(\widetilde{S}, w_2),
\end{eqnarray*}
where $\{i_1', i_2', \ldots,i_{n-2}' \}=[2,n]\backslash \{i\}$.\\

$1\wedge 4$, $f_1=\mathbf{x}_{\sigma}-\mathbf{x}_{\varepsilon}$,
$f_4=\mathbf{x}_{\varepsilon}x_{i_1}x_{i_2}\cdots
x_{i_{m+1}}-\mathbf{x}_{\varepsilon}\underline{x_{i_1}x_{i_2}\cdots
x_{i_{m+1}}}$, $w=x_{j_1}x_{j_2}\cdots
x_{j_t}\mathbf{x}_{\varepsilon}x_{i_1}x_{i_2}\cdots x_{i_{m+1}}$,
$\mathbf{x}_{\sigma}=x_{j_1}x_{j_2}\cdots x_{j_t}x_1x_2\cdots
x_{n-t}$, $\{j_1, j_2,\ldots, j_t\}=[n-t+1, n]$, $2\leq i_1, i_2,
\cdots,i_{m+1}\leq n,$ $m \geq 1$, $1\leq t\leq n-1$,
$x_{i_1}x_{i_2}\cdots x_{i_{m+1}}> \underline{x_{i_1}x_{i_2}\cdots
x_{i_{m+1}}}$.

\begin{eqnarray*}
1\wedge 4 &=&f_1x_{n-t+1}\cdots x_{n}x_{i_1}x_{i_2}\cdots
x_{i_{m+1}}-x_{j_1}x_{j_2}\cdots
x_{j_t}f_4\\
&\equiv&\mathbf{x}_{\varepsilon}x_{n-t+1}\cdots x_{n}x_{i_1}x_{i_2}\cdots
x_{i_{m+1}}-x_{j_1}x_{j_2}\cdots
x_{j_t}\mathbf{x}_{\varepsilon}\underline{x_{i_1}x_{i_2}\cdots
x_{i_{m+1}}}\\
&\equiv&\mathbf{x}_{\varepsilon}\underline{x_{n-t+1}\cdots x_{n}x_{i_1}x_{i_2}\cdots
x_{i_{m+1}}}-\mathbf{x}_{\varepsilon}x_{j_1}x_{j_2}\cdots
x_{j_t}\underline{x_{i_1}x_{i_2}\cdots
x_{i_{m+1}}}\\
&\equiv&\mathbf{x}_{\varepsilon}\underline{x_{n-t+1}\cdots x_{n}x_{i_1}x_{i_2}\cdots
x_{i_{m+1}}}-\mathbf{x}_{\varepsilon}\underline{x_{j_1}x_{j_2}\cdots
x_{j_t}\underline{x_{i_1}x_{i_2}\cdots
x_{i_{m+1}}}}\\
&\equiv&0\ mod(\widetilde{S}, w)
\end{eqnarray*}

$1\wedge 5$, $f_1=\mathbf{x}_{\sigma}-\mathbf{x}_{\varepsilon}$,
$f_5=\mathbf{x}_{\varepsilon}x_{i_1}x_{i_2}\cdots
x_{i_m}x_1-x_1\mathbf{x}_{\varepsilon}x_{i_1}x_{i_2}\cdots x_{i_m}$,
$w=x_{j_1}x_{j_2}\cdots
x_{j_t}\mathbf{x}_{\varepsilon}x_{i_1}x_{i_2}\cdots x_{i_{m}}x_1$,
$\mathbf{x}_{\sigma}=x_{j_1}x_{j_2}\cdots x_{j_t}x_1x_2\cdots
x_{n-t}$,
$\{j_1, j_2,\ldots, j_t\}=[n-t+1, n]$, $2\leq i_1, i_2, \cdots, i_m\leq n$, $m \geq 1$, $1\leq t\leq n-1$.

\begin{eqnarray*}
1\wedge 5 &=&f_1x_{n-t+1}\cdots x_{n}x_{i_1}x_{i_2}\cdots
x_{i_{m}}x_1-x_{j_1}x_{j_2}\cdots
x_{j_t}f_4\\
&\equiv&\mathbf{x}_{\varepsilon}x_{n-t+1}\cdots x_{n}x_{i_1}x_{i_2}\cdots
x_{i_{m}}x_1-x_{j_1}x_{j_2}\cdots
x_{j_t}x_1\mathbf{x}_{\varepsilon}\underline{x_{i_1}x_{i_2}\cdots
x_{i_{m}}}\\
&\equiv&x_1\mathbf{x}_{\varepsilon}\underline{x_{n-t+1}\cdots x_{n}x_{i_1}x_{i_2}\cdots
x_{i_{m}}}-x_1\mathbf{x}_{\varepsilon}x_{j_1}x_{j_2}\cdots
x_{j_t}\underline{x_{i_1}x_{i_2}\cdots
x_{i_{m}}}\\
&\equiv&x_1\mathbf{x}_{\varepsilon}\underline{x_{n-t+1}\cdots x_{n}x_{i_1}x_{i_2}\cdots
x_{i_{m}}}-x_1\mathbf{x}_{\varepsilon}\underline{x_{j_1}x_{j_2}\cdots
x_{j_t}\underline{x_{i_1}x_{i_2}\cdots
x_{i_{m}}}}\\
&\equiv&0\ mod(\widetilde{S}, w)
\end{eqnarray*}

$2\wedge 1$,
$f_2=x_i\mathbf{x}_{\varepsilon}-\mathbf{x}_{\varepsilon}x_i,$
$f_1=\mathbf{x}_{\sigma}-\mathbf{x}_{\varepsilon}$,
$w=x_i\mathbf{x}_{\varepsilon}x_{j_1}x_{j_2}\cdots x_{j_t}$,
$\mathbf{x}_{\sigma}=x_{t+1}\cdots x_n x_{j_1}x_{j_2}\cdots x_{j_t}$,
$\{j_1,j_2, \ldots, j_t \}=\mathbf{n}\backslash [t+1, n]$, $2\leq i\leq n$,  $1\leq t\leq n-1$.
 \begin{eqnarray*}
2\wedge 1&=&f_2x_{j_1}x_{j_2}\cdots
x_{j_t}-x_ix_1x_2\cdots x_{t}f_1\\
&\equiv&\mathbf{x}_{\varepsilon}x_ix_{j_1}x_{j_2}\cdots
x_{j_t}-x_ix_1x_2\cdots x_{t}\mathbf{x}_{\varepsilon}\\
&\equiv&x_1\mathbf{x}_{\varepsilon}x_ix_{j_2'}x_{j_3'}\cdots
x_{j_t'}-x_1\mathbf{x}_{\varepsilon}x_ix_1x_2\cdots x_{t}\\
&\equiv&x_1\mathbf{x}_{\varepsilon}\underline{x_ix_{j_2'}x_{j_3'}\cdots
x_{j_t'}}-x_1\mathbf{x}_{\varepsilon}\underline{x_ix_1x_2\cdots x_{t}}\\
&\equiv&0\ mod(\widetilde{S}, w)
\end{eqnarray*}
where $\{ j_2', j_3',\ldots, j_t' \}=[2,t]$.\\

$2\wedge 2$,
$f_2=x_i\mathbf{x}_{\varepsilon}-\mathbf{x}_{\varepsilon}x_i,$,
$f_2'=x_n\mathbf{x}_{\varepsilon}-\mathbf{x}_{\varepsilon}x_n$,
$w=x_i\mathbf{x}_{\varepsilon}\mathbf{x}_{\varepsilon}$, $2\leq
i\leq  n$.
\begin{eqnarray*}
2\wedge 2&=&f_2\mathbf{x}_{\varepsilon}-x_ix_1\cdots x_{n-1}f_2'\\
&\equiv&\mathbf{x}_{\varepsilon} x_i\mathbf{x}_{\varepsilon}-x_ix_1\cdots x_{n-1}\mathbf{x}_{\varepsilon}x_n\\
&\equiv&x_1\mathbf{x}_{\varepsilon} x_ix_2\cdots x_{n}-x_1\mathbf{x}_{\varepsilon}x_ix_2\cdots x_{n-1}x_n\\
&\equiv&0\ mod(\widetilde{S}, w)
\end{eqnarray*}

$2\wedge 3$,
$f_2=x_i\mathbf{x}_{\varepsilon}-\mathbf{x}_{\varepsilon}x_i$,
$f_3=x_nx_1^m\mathbf{x}_{\varepsilon}-x_1^m\mathbf{x}_{\varepsilon}x_n$,
$w=x_i\mathbf{x}_{\varepsilon}x_1^m\mathbf{x}_{\varepsilon}$,  $2 \leq i\leq n$, $m\geq 1$.
\begin{eqnarray*}
2\wedge 3&=&f_2x_1^m\mathbf{x}_{\varepsilon}-x_ix_1\cdots x_{n-1}f_3\\
&\equiv&\mathbf{x}_{\varepsilon} x_ix_1^m\mathbf{x}_{\varepsilon}-x_ix_1\cdots x_{n-1}x_1^m\mathbf{x}_{\varepsilon}x_n\\
&\equiv&x_1^{m+1}\mathbf{x}_{\varepsilon} x_ix_2\cdots x_n-x_1^{m+1}\mathbf{x}_{\varepsilon} x_ix_2\cdots x_n\\
&\equiv&0\ mod(\widetilde{S}, w)
\end{eqnarray*}

$2\wedge 4$,
$f_2=x_i\mathbf{x}_{\varepsilon}-\mathbf{x}_{\varepsilon}x_i$,
$f_4=\mathbf{x}_{\varepsilon}x_{i_1}x_{i_2}\cdots x_{i_{m+1}}-\mathbf{x}_{\varepsilon}\underline{x_{i_1}x_{i_2}\cdots x_{i_{m+1}}}$,
$w=x_i\mathbf{x}_{\varepsilon}x_{i_1}x_{i_2}\cdots x_{i_{m+1}}$,
$m\geq 1$, $2\leq i, i_1,i_2,\ldots,i_{m+1}\leq n$,
$x_{i_1}x_{i_2}\cdots x_{i_{m+1}}> \underline{x_{i_1}x_{i_2}\cdots
x_{i_{m+1}}}$.

\begin{eqnarray*}
2\wedge 4&=&f_2x_{i_1}x_{i_2}\cdots x_{i_{m+1}}-x_if_4\\
&\equiv&\mathbf{x}_{\varepsilon}x_ix_{i_1}x_{i_2}\cdots x_{i_{m+1}}-x_i\mathbf{x}_{\varepsilon}\underline{x_{i_1}x_{i_2}\cdots x_{i_{m+1}}}\\
&\equiv&\mathbf{x}_{\varepsilon}\underline{x_ix_{i_1}x_{i_2}\cdots x_{i_{m+1}}}-\mathbf{x}_{\varepsilon}x_i\underline{x_{i_1}x_{i_2}\cdots x_{i_{m+1}}}\\
&\equiv&\mathbf{x}_{\varepsilon}\underline{x_ix_{i_1}x_{i_2}\cdots x_{i_{m+1}}}-\mathbf{x}_{\varepsilon}\underline{x_i\underline{x_{i_1}x_{i_2}\cdots x_{i_{m+1}}}}\\
&\equiv&0\ mod(\widetilde{S}, w)
\end{eqnarray*}

$2\wedge 5$,
$f_2=x_i\mathbf{x}_{\varepsilon}-\mathbf{x}_{\varepsilon}x_i$,
$f_5=\mathbf{x}_{\varepsilon}x_{i_1}x_{i_2}\cdots
x_{i_m}x_1-x_1\mathbf{x}_{\varepsilon}x_{i_1}x_{i_2}\cdots x_{i_m}$,
 $w=x_i\mathbf{x}_{\varepsilon}x_{i_1}x_{i_2}\cdots x_{i_m}x_1$, $m\geq 1$,
 $2\leq i, i_1, i_2, \ldots, i_{m}\leq n$.
\begin{eqnarray*}
2\wedge 5&=&f_2x_{i_1}x_{i_2}\cdots x_{i_{m}}x_1-x_if_5\\
&\equiv& \mathbf{x}_{\varepsilon}x_ix_{i_1}x_{i_2}\cdots x_{i_{m}}x_1-x_ix_1\mathbf{x}_{\varepsilon}x_{i_1}x_{i_2}\cdots x_{i_m}\\
&\equiv& \mathbf{x}_{\varepsilon}x_ix_{i_1}x_{i_2}\cdots x_{i_{m}}x_1-x_1 \mathbf{x}_{\varepsilon}x_i x_{i_1}x_{i_2}\cdots x_{i_m}\\
&\equiv& x_1\mathbf{x}_{\varepsilon}x_ix_{i_1}x_{i_2}\cdots x_{i_{m}}-x_1 \mathbf{x}_{\varepsilon}x_i x_{i_1}x_{i_2}\cdots x_{i_m}\\
&\equiv& 0\ mod (\widetilde{S}, w).
\end{eqnarray*}
$3\wedge 1$,
$f_3=x_ix_1^m\mathbf{x}_{\varepsilon}-x_1^m\mathbf{x}_{\varepsilon}x_i$,
 $f_1=\mathbf{x}_{\sigma}-\mathbf{x}_{\varepsilon}$, $w=x_ix_1^m\mathbf{x}_{\varepsilon}x_{i_1}\cdots x_{i_t}$,
$\mathbf{x}_{\sigma}=x_{t+1}x_{t+2}\cdots x_{n}x_{i_1}\cdots x_{i_t}$,
 $\{i_1,i_2,\ldots, i_t\}=[1, t]$, $1\leq t\leq n-1$, $2\leq i \leq n$, $m\geq 1$. \\
 \begin{eqnarray*}
3\wedge 1&=&f_3x_{i_1}\cdots x_{i_t}-x_ix_1^mx_1\cdots x_tf_1\\
&\equiv&x_1^m\mathbf{x}_{\varepsilon}x_ix_{i_1}\cdots x_{i_t}-x_ix_1^mx_1\cdots x_t\mathbf{x}_{\varepsilon}\\
&\equiv&x_1^{m+1}\mathbf{x}_{\varepsilon}\underline{x_ix_{i_2'}\cdots x_{i_t'}}-x_1^{m+1}\mathbf{x}_{\varepsilon}\underline{x_ix_2\cdots x_t}\\
&\equiv& 0\ mod (\widetilde{S}, w).
\end{eqnarray*}
where $\{i_2', i_3', \cdots, i_t'\}=[2,t]$.

$3\wedge 2$,
$f_3=x_ix_1^m\mathbf{x}_{\varepsilon}-x_1^m\mathbf{x}_{\varepsilon}x_i$,
$f_2=x_n\mathbf{x}_{\varepsilon}-\mathbf{x}_{\varepsilon}x_n$,
$w=x_ix_1^m\mathbf{x}_{\varepsilon}\mathbf{x}_{\varepsilon}$, $2\leq
i \leq n$, $m\geq 1$.
 \begin{eqnarray*}
3\wedge 2&=&f_3\mathbf{x}_{\varepsilon}-x_ix_1^mx_1\cdots x_{n-1}f_2\\
&\equiv&x_1^m\mathbf{x}_{\varepsilon}x_i\mathbf{x}_{\varepsilon}-x_ix_1^mx_1\cdots x_{n-1}\mathbf{x}_{\varepsilon}x_n\\
&\equiv&x_1^{m+1}\mathbf{x}_{\varepsilon}x_ix_2\cdots x_n-x_1^{m+1}\mathbf{x}_{\varepsilon}x_ix_2\cdots x_n\\
&\equiv&  0 \ mod (\widetilde{S}, w).
\end{eqnarray*}

$3\wedge 3$,
$f_3=x_ix_1^m\mathbf{x}_{\varepsilon}-x_1^m\mathbf{x}_{\varepsilon}x_i$,
  $f_3'=x_nx_1^{m_1}\mathbf{x}_{\varepsilon}-x_1^{m_1}\mathbf{x}_{\varepsilon}x_n$,
  $w=x_ix_1^{m}\mathbf{x}_{\varepsilon}x_1^{m_1}\mathbf{x}_{\varepsilon}$,
$2\leq i \leq n$, $m, m_1\geq 1$.

 \begin{eqnarray*}
3\wedge 3&=&f_3x_1^{m_1}\mathbf{x}_{\varepsilon}-x_ix_1^mx_1\cdots x_{n-1}f_3'\\
&\equiv&x_1^m\mathbf{x}_{\varepsilon}x_ix_1^{m_1}\mathbf{x}_{\varepsilon}-x_ix_1^mx_1\cdots x_{n-1}x_1^{m_1}\mathbf{x}_{\varepsilon}x_n\\
&\equiv&x_1^{m+m_1+1}\mathbf{x}_{\varepsilon}x_ix_2\cdots x_n-x_1^{m+m_1+1}\mathbf{x}_{\varepsilon}x_ix_2\cdots x_n\\
&\equiv& 0 \ mod (\widetilde{S}, w).
\end{eqnarray*}

$3\wedge 4$,
$f_3=x_ix_1^m\mathbf{x}_{\varepsilon}-x_1^m\mathbf{x}_{\varepsilon}x_i$,
$f_4=\mathbf{x}_{\varepsilon}x_{i_1}\cdots
x_{i_{m_1+1}}-\mathbf{x}_{\varepsilon}\underline{x_{i_1}\cdots
x_{i_{m_1+1}}}$,
$w=x_ix_1^m\mathbf{x}_{\varepsilon}x_{i_1}x_{i_2}\cdots
x_{i_{m_1+1}}$,
 $2\leq i, i_1,i_2,\ldots,i_{m_1+1}\leq n$,\
 $m,m_1\geq 1$, $x_{i_1}\cdots
x_{i_{m_1+1}}>\underline{x_{i_1}\cdots
x_{i_{m_1+1}}}$.

 \begin{eqnarray*}
3\wedge 4&=&f_3x_{i_1}\cdots
x_{i_{m_1+1}}-x_ix_1^mf_4\\
&\equiv&x_1^m\mathbf{x}_{\varepsilon}x_ix_{i_1}\cdots
x_{i_{m_1+1}}-x_ix_1^m\mathbf{x}_{\varepsilon}\underline{x_{i_1}\cdots
x_{i_{m_1+1}}}\\
&\equiv&x_1^m\mathbf{x}_{\varepsilon}\underline{x_ix_{i_1}\cdots
x_{i_{m_1+1}}}-x_1^m\mathbf{x}_{\varepsilon}\underline{x_i x_{i_1}\cdots
x_{i_{m_1+1}}}\\
&\equiv& 0\  mod (\widetilde{S}, w).
\end{eqnarray*}

$3\wedge 5$,
$f_3=x_ix_1^m\mathbf{x}_{\varepsilon}-x_1^m\mathbf{x}_{\varepsilon}x_i$,
 $f_5=\mathbf{x}_{\varepsilon}x_{i_1}x_{i_2}\cdots x_{i_{m_1}}x_1-x_1\mathbf{x}_{\varepsilon}x_{i_1}x_{i_2}\cdots x_{i_{m_1}}$,
$w=x_ix_1^m\mathbf{x}_{\varepsilon}x_{i_1}x_{i_2}\cdots x_{i_{m_1}}x_1$,
 $2\leq i, i_1,i_2,\ldots,i_{m_1}\leq n$,   $m,m_1\geq 1$.\\

\begin{eqnarray*}
3\wedge 5&=&f_3x_{i_1}x_{i_2}\cdots x_{i_{m_1}}x_1-x_ix_1^mf_5\\
&\equiv&x_1^m\mathbf{x}_{\varepsilon}x_ix_{i_1}x_{i_2}\cdots x_{i_{m_1}}x_1-x_ix_1^mx_1\mathbf{x}_{\varepsilon}\underline{x_{i_1}\cdots
x_{i_{m_1}}}\\
&\equiv&x_1^{m+1}\mathbf{x}_{\varepsilon}\underline{x_ix_{i_1}\cdots
x_{i_{m_1}}}-x_1^{m+1}\mathbf{x}_{\varepsilon}\underline{x_i x_{i_1}\cdots
x_{i_{m_1}}}\\
&\equiv& 0 \ mod (\widetilde{S}, w).
\end{eqnarray*}
$4\wedge 1$, there are two cases. Let
$f_4=\mathbf{x}_{\varepsilon}x_{i_1}x_{i_2}\cdots
x_{i_{m+1}}-\mathbf{x}_{\varepsilon}\underline{x_{i_1}x_{i_2}\cdots
x_{i_{m+1}}}$, $f_1=\mathbf{x}_{\sigma}-\mathbf{x}_{\varepsilon}$,
  $m\geq 1$, $2\leq i_1, i_2,\cdots,i_{m+1}\leq n$,
  $x_{i_1}x_{i_2}\cdots x_{i_{m+1}}>\underline{x_{i_1}x_{i_2}\cdots
x_{i_{m+1}}}$.

$w_1=\mathbf{x}_{\varepsilon}x_{i_1}\cdots
x_{i_{m+1}}x_{j_1}x_{j_2}\cdots x_{j_{t}}$,
$\mathbf{x}_{\sigma}=x_{i_{m+2-n+t}} \cdots
x_{i_{m+1}}x_{j_1}x_{j_2}\cdots x_{j_{t}}$,\\
 $\{ j_1, j_2, \ldots,
j_{t}\}=\mathbf{n}\backslash \{ i_{m+2-n+t}, \ldots, i_{m+1}\}$,

$w_2=\mathbf{x}_{\varepsilon}x_{i_1} \cdots x_{i_{m+1}}x_{j_1}\cdots
x_{j_{t-m-1}}$, $\mathbf{x}_{\sigma}=x_{t+1}\cdots x_n x_{i_1}
\cdots x_{i_{m+1}}x_{j_1}\cdots x_{j_{t-m-1}}$, $2\leq i_1,
i_2,\cdots,i_{m+1}\leq n$, $\{ j_1, j_2, \ldots,
j_{t-m-1}\}=\mathbf{n}\backslash ([t+1, n]\cup \{ i_{1},
i_{2},\ldots, i_{m+1}\})$,
 $t-m-1\geq 1$.

\begin{eqnarray*}
4\wedge 1&=&f_4x_{j_1}x_{j_2}\cdots x_{j_{t}}-\mathbf{x}_{\varepsilon}x_{i_1}\cdots x_{i_{m+1-n+t}}f_1\\
&\equiv&\mathbf{x}_{\varepsilon}\underline{x_{i_1}x_{i_2}\cdots
x_{i_{m+1}}}x_{j_1}x_{j_2}\cdots x_{j_{t}}-\mathbf{x}_{\varepsilon}x_{i_1}\cdots x_{i_{m+1-n+t}}\mathbf{x}_{\varepsilon}\\
&\equiv& x_1\mathbf{x}_{\varepsilon}\underline{x_{i_1}x_{i_2}\cdots
x_{i_{m+1}}x_{j_2'}x_{j_3'}\cdots x_{j_{t'}}}-x_1\mathbf{x}_{\varepsilon}\underline{x_{i_1}\cdots x_{i_{m+1-n+t}}x_2\cdots x_n}\\
&\equiv& 0 \ mod (\widetilde{S}, w_1).
\end{eqnarray*}
where $\{j_2',j_3',\ldots,j_{t'}\}=\mathbf{n}\backslash \{ i_{m+2-n+t}, \ldots, i_{m+1}, 1\}$.

\begin{eqnarray*}
4\wedge 1&=&f_4x_{j_1}x_{j_2}\cdots x_{j_{t-m-1}}-x_{1}\cdots x_{t}f_1\\
&\equiv&\mathbf{x}_{\varepsilon}\underline{x_{i_1}x_{i_2}\cdots
x_{i_{m+1}}}x_{j_1}x_{j_2}\cdots x_{j_{t-m-1}}-x_{1}\cdots x_{t}\mathbf{x}_{\varepsilon}\\
&\equiv&x_{1}\mathbf{x}_{\varepsilon}\underline{x_{i_1}x_{i_2}\cdots
x_{i_{m+1}}x_{j_2'}x_{j_3'}\cdots x_{j_{t-m-1}'}}-x_{1}\mathbf{x}_{\varepsilon}x_{2}\cdots x_{t}\\
&\equiv& 0 \ mod (\widetilde{S}, w_1).
\end{eqnarray*}
where $\{j_2',j_3',\ldots,j_{t-m-1}'\}=[2,t]\backslash\{ i_{1},
i_{2},\ldots, i_{m+1}\}$.

$4\wedge 2$, $f_4=\mathbf{x}_{\varepsilon}x_{i_1}x_{i_2}\cdots
x_{i_{m+1}}-\mathbf{x}_{\varepsilon}\underline{x_{i_1}x_{i_2}\cdots
x_{i_{m+1}}}$,
$f_2=x_{i_{m+1}}\mathbf{x}_{\varepsilon}-\mathbf{x}_{\varepsilon}x_{i_{m+1}}$,
$w=\mathbf{x}_{\varepsilon}x_{i_1}x_{i_2}\cdots
x_{i_{m+1}}\mathbf{x}_{\varepsilon}$, $m\geq 1$, $2\leq
i_1,\cdots,i_{m+1}\leq n$, $x_{i_1}x_{i_2}\cdots
x_{i_{m+1}}>\underline{x_{i_1}x_{i_2}\cdots x_{i_{m+1}}}$.

\begin{eqnarray*}
4\wedge 2&=&f_4\mathbf{x}_{\varepsilon}-\mathbf{x}_{\varepsilon}x_{i_1}x_{i_2}\cdots
x_{i_{m}}f_2\\
&\equiv&\mathbf{x}_{\varepsilon}\underline{x_{i_1}x_{i_2}\cdots
x_{i_{m+1}}}\mathbf{x}_{\varepsilon}-\mathbf{x}_{\varepsilon}x_{i_1}x_{i_2}\cdots
x_{i_{m}}\mathbf{x}_{\varepsilon}x_{i_{m+1}}\\
&\equiv& x_1\mathbf{x}_{\varepsilon}\underline{x_{i_1}x_{i_2}\cdots
x_{i_{m+1}}x_2\cdots x_n}-x_1\mathbf{x}_{\varepsilon}\underline{x_{i_1}x_{i_2}\cdots
x_{i_{m}}x_2\cdots x_nx_{i_{m+1}}}\\
&\equiv& 0 \ mod (\widetilde{S}, w).
\end{eqnarray*}

$4\wedge 3$, $f_4=\mathbf{x}_{\varepsilon}x_{i_1}x_{i_2}\cdots
x_{i_{m+1}}-\mathbf{x}_{\varepsilon}\underline{x_{i_1}x_{i_2}\cdots
x_{i_{m+1}}}$,
$f_3=x_{i_{m+1}}x_1^{m_1}\mathbf{x}_{\varepsilon}-x_1^{m_1}\mathbf{x}_{\varepsilon}x_{i_{m+1}}$,
$w=\mathbf{x}_{\varepsilon}x_{i_1}x_{i_2}\cdots
x_{i_{m+1}}x_1^{m_1}\mathbf{x}_{\varepsilon}$, $2\leq
i_1,\cdots,i_{m+1}\leq n$, $m_1, m\geq 1$, $x_{i_1}x_{i_2}\cdots
x_{i_{m+1}}>\underline{x_{i_1}x_{i_2}\cdots x_{i_{m+1}}}$.

\begin{eqnarray*}
4\wedge 3&=&f_4x_1^m\mathbf{x}_{\varepsilon}-\mathbf{x}_{\varepsilon}x_{i_1}x_{i_2}\cdots
x_{i_{m}}f_3\\
&\equiv&\mathbf{x}_{\varepsilon}\underline{x_{i_1}x_{i_2}\cdots
x_{i_{m+1}}}x_1^m\mathbf{x}_{\varepsilon}-\mathbf{x}_{\varepsilon}x_{i_1}x_{i_2}\cdots
x_{i_{m}}x_1^m\mathbf{x}_{\varepsilon}x_{i_{m+1}}\\
&\equiv&x_1^{m+1}\mathbf{x}_{\varepsilon}\underline{x_{i_1}\cdots
x_{i_{m+1}}}x_2\cdots x_n-x_1^{m+1}\mathbf{x}_{\varepsilon}x_{i_1}\cdots
x_{i_{m}}x_2\cdots x_nx_{i_{m+1}}\\
&\equiv&x_1^{m+1}\mathbf{x}_{\varepsilon}\underline{x_{i_1}\cdots
x_{i_{m+1}}x_2\cdots x_n}-x_1^{m+1}\mathbf{x}_{\varepsilon}\underline{x_{i_1}\cdots
x_{i_{m}}x_2\cdots x_nx_{i_{m+1}}}\\
&\equiv& 0 \ mod (\widetilde{S}, w).
\end{eqnarray*}

$5\wedge 1$,  there are two cases. Let
 $f_5=\mathbf{x}_{\varepsilon}x_{i_1}x_{i_2}\cdots
x_{i_m}x_1-x_1\mathbf{x}_{\varepsilon}x_{i_1}x_{i_2}\cdots x_{i_m}$,
$f_1=\mathbf{x}_{\sigma}-\mathbf{x}_{\varepsilon}$, $2\leq i_1, i_2,
\cdots , i_m\leq n$, $m\geq 1$.

$w_1=\mathbf{x}_{\varepsilon}x_{i_1}\cdots x_{i_m}x_1x_{j_1}\cdots
x_{j_{n+t-m-2}}$, $\mathbf{x}_{\sigma}=x_{i_{t}}x_{i_{t+1}}\cdots
x_{i_m}x_1x_{j_1}\cdots x_{j_{n+t-m-2}}$, $2\leq i_1, i_2, \cdots ,
i_m\leq n$, $\{ j_1,j_2, \ldots, j_{n+t-m-2} \}
=\mathbf{n}\backslash(\{i_{t},i_{t+1},\ldots,i_m\}\cup \{1\})$.

$w_2=\mathbf{x}_{\varepsilon}x_{i_1}\cdots x_{i_m}x_1$,
$\mathbf{x}_{\sigma}=x_{t+1}\cdots x_nx_{i_1}x_{i_2}\cdots
x_{i_m}x_1x_{j_1}\cdots x_{j_{t-m}}$, \\  $2\leq i_1, i_2, \cdots,
i_m\leq n$, $\{ j_1,j_2, \ldots, j_{t-m}
\}=\mathbf{n}\backslash(\{i_{1},i_{2},\ldots,i_m\}\cup [t+1, n]\cup
\{1\})$, $ 1\leq m \leq n-2$, $t-m\geq 0$.

\begin{eqnarray*}
5\wedge 1&=&f_5x_{j_1}\cdots
x_{j_{n+t-m-2}}-\mathbf{x}_{\varepsilon}x_{i_1}x_{i_2}\cdots
x_{i_{t-1}}f_1\\
&\equiv&x_1\mathbf{x}_{\varepsilon}x_{i_1}x_{i_2}\cdots x_{i_m}x_{j_1}\cdots
x_{j_{n+t-m-2}}-\mathbf{x}_{\varepsilon}x_{i_1}\cdots
x_{i_{t-1}}\mathbf{x}_{\varepsilon}\\
&\equiv&x_1\mathbf{x}_{\varepsilon}x_{i_1}x_{i_2}\cdots x_{i_m}x_{j_1}\cdots
x_{j_{n+t-m-2}}-x_1\mathbf{x}_{\varepsilon}x_{i_1}\cdots
x_{i_{t-1}}x_2\cdots x_n\\
&\equiv&x_1\mathbf{x}_{\varepsilon}\underline{x_{i_1}x_{i_2}\cdots x_{i_m}x_{j_1}\cdots
x_{j_{n+t-m-2}}}-x_1\mathbf{x}_{\varepsilon}\underline{x_{i_1}\cdots
x_{i_{t-1}}x_2\cdots x_n}\\
&\equiv& 0 \ mod (\widetilde{S}, w_1).
\end{eqnarray*}
\begin{eqnarray*}
5\wedge 1&=&f_5x_{j_1}\cdots
x_{j_{t-m}}-x_{1}x_{2}\cdots
x_{t}f_1\\
&\equiv&x_1\mathbf{x}_{\varepsilon}x_{i_1}x_{i_2}\cdots x_{i_m}x_{j_1}\cdots
x_{j_{t-m}}-x_{1}x_{2}\cdots
x_{t}\mathbf{x}_{\varepsilon}\\
&\equiv&x_1\mathbf{x}_{\varepsilon}x_{i_1}x_{i_2}\cdots x_{i_m}x_{j_1}\cdots
x_{j_{t-m}}-x_{1}\mathbf{x}_{\varepsilon}x_{2}\cdots
x_{t}\\
&\equiv&x_1\mathbf{x}_{\varepsilon}\underline{x_{i_1}x_{i_2}\cdots x_{i_m}x_{j_1}\cdots
x_{j_{t-m}}}-x_{1}\mathbf{x}_{\varepsilon}\underline{x_{2}\cdots
x_{t}}\\
&\equiv& 0 \ mod (\widetilde{S}, w_2).
\end{eqnarray*}

$5\wedge 2$,  $f_5=\mathbf{x}_{\varepsilon}x_{i_1}x_{i_2}\cdots
x_{i_m}x_1-x_1\mathbf{x}_{\varepsilon}x_{i_1}x_{i_2}\cdots x_{i_m}$,
$f_2=x_{i_m}\mathbf{x}_{\varepsilon}-\mathbf{x}_{\varepsilon}x_{i_m}$
$w=\mathbf{x}_{\varepsilon}x_{i_1}x_{i_2}\cdots
x_{i_m}\mathbf{x}_{\varepsilon}$,
$2\leq i_1, i_2, \cdots, i_{m}\leq n$, $m\geq 1$. \\

\begin{eqnarray*}
5\wedge 2&=&f_5x_{2}\cdots
x_{n}-\mathbf{x}_{\varepsilon}x_{i_1}\cdots
x_{i_{m-1}}f_2\\
&\equiv&x_1\mathbf{x}_{\varepsilon}x_{i_1}x_{i_2}\cdots x_{i_m}x_{2}\cdots
x_{n}-\mathbf{x}_{\varepsilon}x_{i_1}\cdots
x_{i_{m-1}}\mathbf{x}_{\varepsilon}x_{i_m}\\
&\equiv&x_1\mathbf{x}_{\varepsilon}x_{i_1}x_{i_2}\cdots x_{i_m}x_{2}\cdots
x_{n}-x_1\mathbf{x}_{\varepsilon}x_{i_1}\cdots
x_{i_{m-1}}x_2\cdots x_nx_{i_m}\\
&\equiv&x_1\mathbf{x}_{\varepsilon}\underline{x_{i_1}x_{i_2}\cdots x_{i_m}x_{2}\cdots
x_{n}}-x_1\mathbf{x}_{\varepsilon}\underline{x_{i_1}\cdots
x_{i_{m-1}}x_2\cdots x_nx_{i_m}}\\
&\equiv& 0 \ mod (\widetilde{S}, w).
\end{eqnarray*}
$5\wedge 3$, $f_5=\mathbf{x}_{\varepsilon}x_{i_1}x_{i_2}\cdots
x_{i_m}x_1-x_1\mathbf{x}_{\varepsilon}x_{i_1}x_{i_2}\cdots x_{i_m}$,
$f_3=
x_{i_m}x_1^{m_1}\mathbf{x}_{\varepsilon}-x_1^{m_1}\mathbf{x}_{\varepsilon}x_{i_m}$,
$w=\mathbf{x}_{\varepsilon}x_{i_1}x_{i_2}\cdots
x_{i_m}x_1^{m_1}\mathbf{x}_{\varepsilon}$,
$2\leq i_1, i_2, \cdots,  i_{m}\leq n$, $m, m_1\geq 1$.

\begin{eqnarray*}
5\wedge 3&=&f_5x_{1}^{m_1-1}\mathbf{x}_{\varepsilon}-\mathbf{x}_{\varepsilon}x_{i_1}\cdots
x_{i_{m-1}}f_3\\
&\equiv&x_1\mathbf{x}_{\varepsilon}x_{i_1}x_{i_2}\cdots x_{i_m}x_{1}^{m_1-1}\mathbf{x}_{\varepsilon}-\mathbf{x}_{\varepsilon}x_{i_1}\cdots
x_{i_{m-1}}x_1^{m_1}\mathbf{x}_{\varepsilon}\\
&\equiv&x_{1}^{m_1+1}\mathbf{x}_{\varepsilon}x_{i_1}\cdots x_{i_m}x_2\cdots  x_n -x_{1}^{m_1+1}\mathbf{x}_{\varepsilon}x_{i_1}\cdots x_{i_m}x_2\cdots  x_n\\
&\equiv& 0 \ mod (\widetilde{S}, w).
\end{eqnarray*}
$5\wedge 4$,
 $f_5=\mathbf{x}_{\varepsilon}x_{i_1}x_{i_2}\cdots
x_{i_m}x_1-x_1\mathbf{x}_{\varepsilon}x_{i_1}x_{i_2}\cdots x_{i_m}$,
$f_4=\mathbf{x}_{\varepsilon}x_{j_1}x_{j_2}\cdots x_{j_{m_1+1}}$ \\
$-\mathbf{x}_{\varepsilon}\underline{x_{j_1}x_{j_2}\cdots
x_{j_{m_1+1}}}$, $w=\mathbf{x}_{\varepsilon}x_{i_1}x_{i_2}\cdots
x_{i_m}\mathbf{x}_{\varepsilon}x_{j_1}x_{j_2}\cdots x_{j_{m_1+1}}$,
$2\leq i_1, i_2,  \cdots,  i_m,\\ j_1,j_2,\cdots, j_{m_1+1}\leq n,$
 $ m, m_1\geq 1$, $x_{j_1}x_{j_2}\cdots x_{j_{m_1+1}}>\underline{x_{j_1}x_{j_2}\cdots x_{j_{m_1+1}}}$.

 \begin{eqnarray*}
 &&5\wedge 4\\
&=&f_5x_{2}\cdots x_nx_{j_1}x_{j_2}\cdots x_{j_{m_1+1}}-\mathbf{x}_{\varepsilon}x_{i_1}\cdots
x_{i_{m}}f_4\\
&\equiv&x_1\mathbf{x}_{\varepsilon}x_{i_1}\cdots x_{i_m}x_{2}\cdots x_nx_{j_1}\cdots x_{j_{m_1+1}}-\mathbf{x}_{\varepsilon}x_{i_1}\cdots
x_{i_{m}}\mathbf{x}_{\varepsilon}x_{j_1}\cdots x_{j_{m_1+1}}\\
&\equiv&x_1\mathbf{x}_{\varepsilon}x_{i_1}\cdots x_{i_m}x_{2}\cdots x_nx_{j_1}\cdots x_{j_{m_1+1}}-x_1\mathbf{x}_{\varepsilon}x_{i_1}\cdots x_{i_m}x_{2}\cdots x_nx_{j_1}\cdots x_{j_{m_1+1}}\\
&\equiv& 0 \ mod (\widetilde{S}, w).
\end{eqnarray*}

$5\wedge 5$, $f_5=\mathbf{x}_{\varepsilon}x_{i_1}x_{i_2}\cdots
x_{i_m}x_1-x_1\mathbf{x}_{\varepsilon}x_{i_1}x_{i_2}\cdots x_{i_m}$,
$f_5'=\mathbf{x}_{\varepsilon}x_{i_1'}x_{i_2'}\cdots x_{i_{m_1}'}x_1$ \\
$-x_1\mathbf{x}_{\varepsilon}x_{i_1'}x_{i_2'}\cdots x_{i_{m_1}'}$,
$w=\mathbf{x}_{\varepsilon}x_{i_1}x_{i_2}\cdots x_{i_m}\mathbf{x}_{\varepsilon}x_{i_1'}x_{i_2'}\cdots x_{i_{m_1}'}x_1$,
$2\leq i_1, i_2,  \cdots,  i_m,\\ i_1', i_2', \cdots,  i_{m_1}'\leq n$.
\begin{eqnarray*}
&&5\wedge 5\\
&=&f_5x_{2}\cdots x_nx_{i_1'}x_{i_2'}\cdots x_{i_{m_1}'}x_1-\mathbf{x}_{\varepsilon}x_{i_1}x_{i_2}\cdots
x_{i_m}f_5'\\
&\equiv&x_1\mathbf{x}_{\varepsilon}x_{i_1}\cdots x_{i_m}x_{2}\cdots x_nx_{i_1'}\cdots x_{i_{m_1}'}x_1
-\mathbf{x}_{\varepsilon}x_{i_1}\cdots
x_{i_m}x_1\mathbf{x}_{\varepsilon}x_{i_1'}\cdots x_{i_{m_1}'}\\
&\equiv&x_1^2\mathbf{x}_{\varepsilon}x_{i_1}\cdots x_{i_m}x_{2}\cdots x_nx_{i_1'}\cdots x_{i_{m_1}'}
-x_1^2\mathbf{x}_{\varepsilon}x_{i_1}\cdots x_{i_m}x_{2}\cdots x_nx_{i_1'}\cdots x_{i_{m_1}'}
\\
&\equiv& 0 \ mod (\widetilde{S}, w).
\end{eqnarray*}

\end{document}